\documentclass[reqno,12pt,letterpaper]{amsart}
\usepackage{amsmath,amssymb,amsthm,graphicx,mathrsfs,url}
\usepackage[usenames,dvipsnames]{color}
\usepackage[colorlinks=true,linkcolor=Red,citecolor=Green]{hyperref}
\usepackage{amsxtra}
\usepackage{wasysym} % for "\varhexagon" macro
\usepackage{graphicx}% for "\scalebox" macro
\usepackage{subcaption}

\usepackage{graphicx,color}

\def\arXiv#1{\href{http://arxiv.org/abs/#1}{arXiv:#1}}
\usepackage{array}
\newcolumntype{P}[1]{>{\centering\arraybackslash}m{#1}}

\def\wrtext#1{\relax\ifmmode{\leavevmode\hbox{#1}}\else{#1}\fi}

\def\?[#1]{\textbf{[#1]}\marginpar{\Large{\textbf{??}}}}

\let\epsilon=\varepsilon % sorry Knuth

\numberwithin{equation}{section}

\DeclareMathOperator{\Spec}{Spec}

\let\Im=\Imag

\DeclareMathOperator{\Op}{Op}

\let\Re=\Real

\DeclareMathOperator{\supp}{supp}

\usepackage{scalerel}

\newcommand\reallywidehat[1]{\arraycolsep=0pt\relax%
\begin{array}{c}
\stretchto{
  \scaleto{
    \scalerel*[\widthof{\ensuremath{#1}}]{\kern-.5pt\bigwedge\kern-.5pt}
    {\rule[-\textheight/2]{1ex}{\textheight}} %WIDTH-LIMITED BIG WEDGE
  }{\textheight} %
}{0.5ex}\\           % THIS SQUEEZES THE WEDGE TO 0.5ex HEIGHT
#1\\                 % THIS STACKS THE WEDGE ATOP THE ARGUMENT
\rule{-1ex}{0ex}
\end{array}
}

\begin{document}

\title[Dissipation for contact Anosov flows]{Optimal enhanced dissipation for contact Anosov flows}

\author{Zhongkai Tao}
\email{ztao@math.berkeley.edu}
\address{Department of Mathematics, University of California,
Berkeley, CA 94720, USA.}

\author{Maciej Zworski}
\email{zworski@math.berkeley.edu}
\address{Department of Mathematics, University of California,
Berkeley, CA 94720, USA.}

\begin{abstract}
We show that for a contact Anosov flow on a compact manifold $ M $, the solutions to 
$ \partial_t u + X u =  \nu \Delta u$, $ \nu > 0 $, where $ X $ is the generator of the flow and 
$ \Delta $, a (negative) Laplacian for some Riemannian metric on $ M $, satisfy
\[    \| u ( t )   -  \underline u \|_{L^2 ( M) } \leq C \nu^{-K} e^{ - \beta  t } \| u( 0 )  \|_{L^2 ( M) }, \]
where $ \underline u $ is the (conserved) average of $ u (0) $ with respect to the contact volume form, and $K$, $\beta$ are fixed positive constants. Since our class of flows includes geodesic flows on 
manifolds of negative curvature, this provides many examples of very precise {\em optimal enhanced dissipation} in the sense of \cite{BBP} and \cite{enhance}. The proof is based on results
about stochastic stability of Pollicott--Ruelle resonances \cite{damped}.
\end{abstract}

\maketitle

\section{Introduction}
Let $M$ be a compact contact manifold such that the Reeb flow has the Anosov property -- 
see \cite[\S 9.1]{nozw} for a brief review. A large class of examples is given by geodesic
flows on compact Riemannian manifolds with negative curvature: $M $ is then given by the co-sphere bundle of the manifold and the Reeb vector field is the generator of the geodesic flow.
\begin{figure}
\includegraphics[width=8cm]{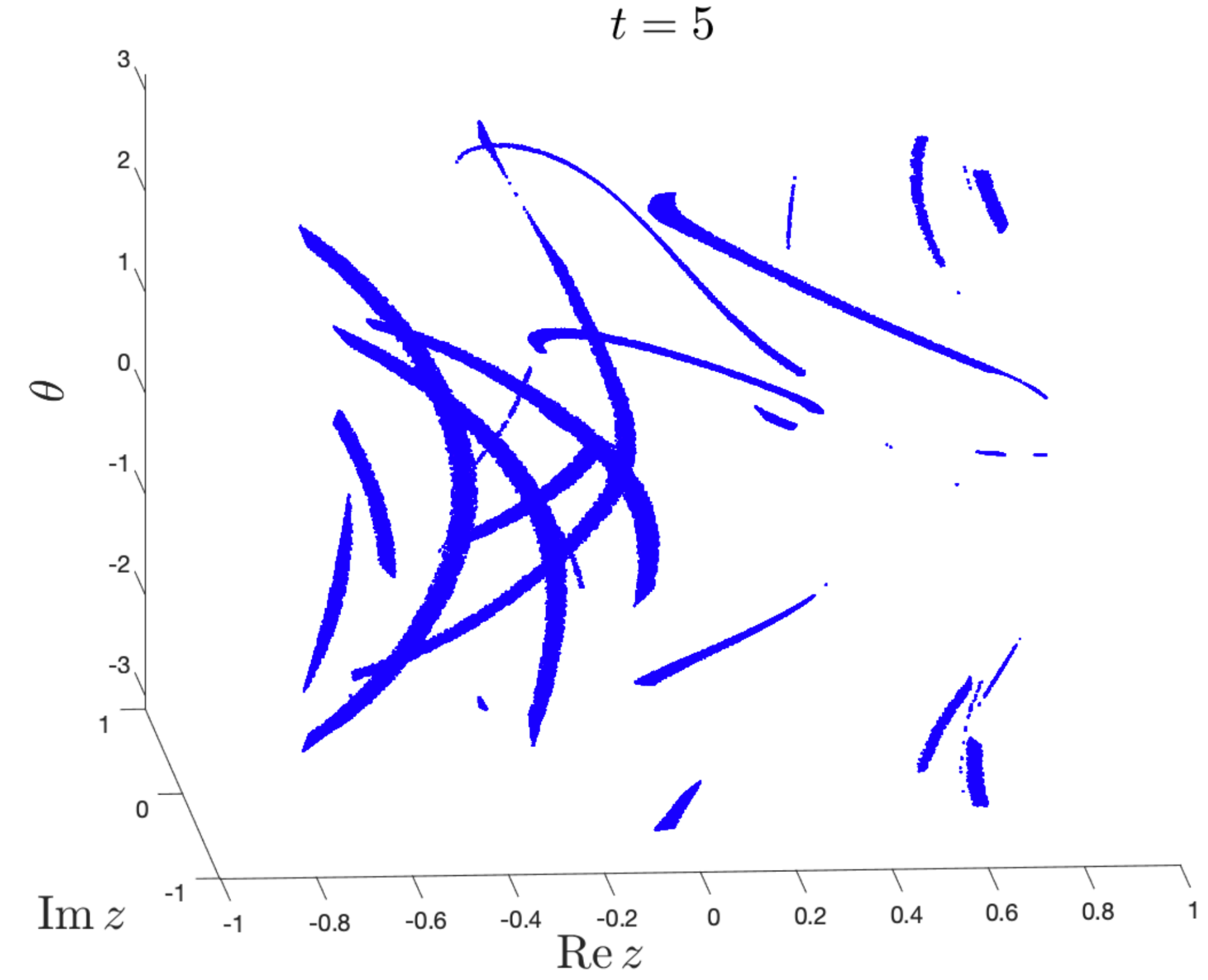}\includegraphics[width=8cm]{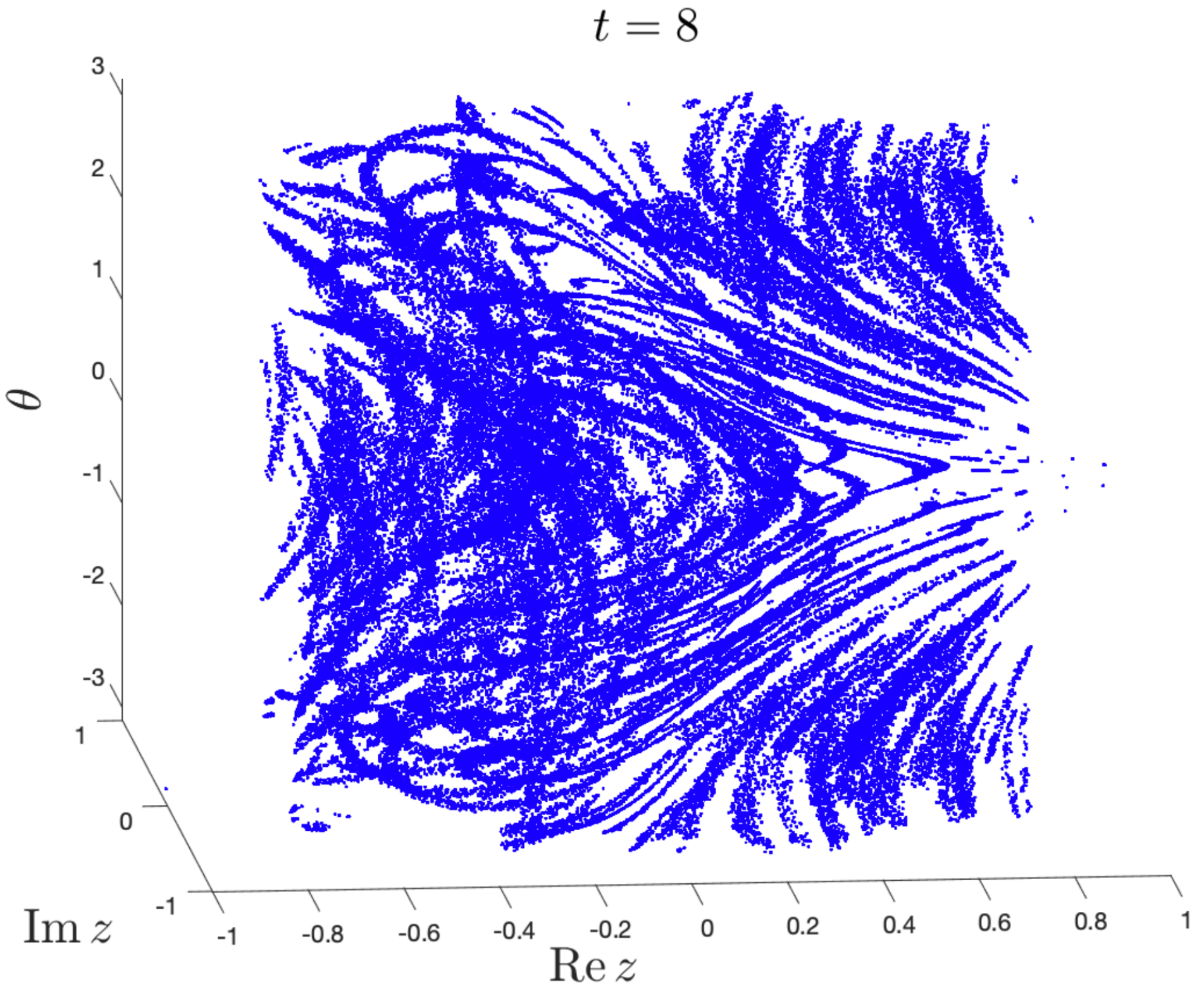}
\caption{\label{f:1} 
An illustration of an exponentially mixing contact flow: the geodesic flow 
on $ S^* \Sigma $ where $  \Sigma = \Gamma \backslash D ( 0, 1 ) $   is the Bolza surface, a genus two 
surface of constant negative curvature; $ z \in D ( 0 , 1 ) $ (the Poincar\'e disc) is the variable in the fundamental domain of $ \Gamma $ and $ e^{ i \theta }  $ is the direction of the geodesic. 
The figures show the evolution of a neighbourhood of $ z = 0 $, $ \theta = \pi/5 $ (approximated by 
$ 10^5$ uniformly distributed points) 
under the flow at times $ t = 5, 8$. The flow is periodic  $ {}\! \! \! \! \mod  \! \Gamma $ in 
$  z $ and ${} \!\!\ \!\! \! \! \mod \!2 \pi$ in $ \theta$ -- see \url{https://math.berkeley.edu/~zworski/bolzamix.mp4} for an animation.  For the flow in the base
see \url{https://rb.gy/xfssmg}. We are grateful to Semyon Dyatlov for help in producing these figures and movies.}
\end{figure}

These flows are exponentially mixing as was shown by Liverani \cite{liver}, who extended 
the work of Dolgopyat \cite{Dol}, with subsequent improvements and generalizations by 
Tsujii \cite{Ts} and Nonnenmacher--Zworski \cite{nozw}. A rough, but relevant here, version of the mixing statement goes as follows: there exist  $ s_0 > 0 $, $ \beta_1 > 0 $, such that for $ f , g \in H^{s_0} ( M ) $ (here
$ H^s $ denotes a Sobolev space) 
\begin{equation}
\label{eq:mixing}
\int_M f ( \rho ) g ( \exp ( t X ) \rho ) dm ( \rho ) - \int_M f \int_M g = e^{ - \beta_1 t } \mathcal O \left( \| f \|_{H^{s_0} } \| g \|_{ H^{s_0 } }  \right)   ,  \end{equation}
where $ d m $ is the contact volume form normalized so that $ \int_M dm = 1 $.

Motivated by recent work of Elgindi--Liss--Mattingly \cite{enhance} we consider the corresponding
convection diffusion equation with $ L^2 $ initial data: let $ \Delta \leq 0 $ be the Laplace operator
for some Riemannian  metric on $M $, and 
\begin{equation}
\label{eq:condi}
\partial_t u + X u = \nu \Delta u , \ \ \  u ( 0 ) \in L^2 ( M ) , \ \ \nu > 0 . 
\end{equation}
Our main result is 

\noindent
{\bf Theorem.} {\em Let $\beta_1$ be the constant in \eqref{eq:condi}, $\gamma_0>0$ be the averaged Lyapunov exponent defined in \eqref{e:lyapunov}, and $0<\beta<\min\{\beta_1,\gamma_0/2\}$. Then there exist constants $ C,  K $ and $\nu_0>0$ %(depending only on $ M$, $ \Delta $ and $ X $)
such that solutions to \eqref{eq:condi} satisfy 
\begin{equation}
\label{eq:theo} \| u ( t )   -  \underline u \|_{L^2 ( M) } \leq C \nu^{-K} e^{ - \beta t } \| u( 0 )  \|_{L^2 ( M) }, \quad t>0, \ 0<\nu<\nu_0
\end{equation}
where $ \underline u $ is the average of $ u $.}

The estimate \eqref{eq:theo} implies a weaker estimate involving only exponential decay
\begin{equation}
\label{eq:weaker}
\| u ( t )   -  \underline u \|_{L^2 ( M) } \leq C  e^{ -  t / \log ( 1/\nu )  } \| u( 0 )  \|_{L^2 ( M) }, \quad t>0, \end{equation}
which is referred to in \cite{enhance} as {\em optimal enhanced dissipation}. In fact,  \eqref{eq:weaker} is
immediate from \eqref{eq:theo} for $t>C\log(1/\nu)$ for some large $C>0$. For $0<t\leq C\log(1/\nu) $ we have (since $ \Re \langle X u , u \rangle = 0 $)
\begin{align*}
    \frac{1}{2}\frac{d}{dt}\|u(t)\|_{L^2}^2=\Re\langle\partial_t u,u\rangle + \Re\langle Xu,u\rangle=\nu\langle \Re\Delta u,u\rangle\leq C\nu \|u\|_{L^2}^2.
\end{align*}
The last step inequality is a consequence of the $L^2$ spectrum of $\Re\Delta=\frac{1}{2}(\Delta+\Delta^*)$ being bounded from above.
Hence,  for $0<t\leq C\log(1/\nu)$,
\begin{equation}
\label{eq:poinc} 
    \|u(t)\|_{L^2}\leq e^{C\nu t}\|u(0)\|_{L^2}\leq C\|u(0)\|_{L^2}.
\end{equation}
We also remark that we could always rescale the Riemannian metric to make its density equal to 
the contact density, in which case $ \Delta $ is self-adjoint for the invariant measure. For
$ u (0 ) $ with $ \underline u = 0 $, that would 
give a bound with $ e^{ - \nu t/ C}$ replacing $ e^{ C \nu t } $ in \eqref{eq:poinc}.

In 
\cite{enhance}  \eqref{eq:weaker} was proved for a specific 
time dependent Lipschitz vector field on the 2-torus designed to induce strong mixing and the authors
stated that ``the only flows we are aware of that achieve this rate are those constructed in
\cite{BBP} from generic solutions to various stochastically forced fluid models." Motivation in 
\cite{enhance} and \cite{BBP} (by Bedrossian--Blumenthal--Punshon-Smith) came from fluid mechanics. Earlier, 
Constantin, Kiselev, Ryzhik and Zlatoš \cite{CKRZ} studied the relationship between the mixing property of $X$ and the qualitative behavior of the enhanced dissipation. A quantitative version was 
obtained by Coti Zelati--Delgadino--Elgindi \cite{CDE}, who achieved a decay rate of $e^{-t/(\log\nu)^2}$ for exponentially mixing flows. We refer to these papers for pointers to the substantial literature on the subject. 

In this paper we consider geometric flows generated by time independent vector fields which
allows a spectral theory point of view. That means considering, in the notation above,
\begin{equation}
\label{eq:defP} 
P_\nu = X - \nu \Delta .
\end{equation}
On $ L ^2 ( M ) $ (with the measure given by the contact volume form) $ i X $ is 
self-adjoint and its spectrum (with the domain $\{u\in L^2: Xu\in L^2\} $) is real.
For $ \nu > 0 $, the operator $ P_\nu $ is elliptic and, with the domain given by 
$ H^2 ( M )$,  its spectrum is discrete. Dyatlov--Zworski
\cite{damped} showed that if the flow generated by $ X $ has the Anosov property 
(see that paper for a definition) then the spectrum of $ P_\nu $ converges, uniformly on 
compact sets, to a discrete set of {\em Pollicott--Ruelle} resonances. The methods of \cite{damped}
are based on the microlocal approach to Anosov flows introduced by Faure--Sj\"ostrand \cite{FS} and
developed in \cite{zeta} and many other works. The Pollicott--Ruelle resonances were introduced by 
Pollicott \cite{poli} and Ruelle \cite{rue}
to describe the power spectrum of correlations. They can be understood as spectrum of $ X $
acting of certain modified spaces -- see \cite{liver} for a dynamical system approach to that and 
references. Estimate \eqref{eq:mixing} for contact flows
 follows from a spectral gap but since the spaces are different
it does not hold with $ L^2 $ bounds. To obtain \eqref{eq:theo} we use the uniformity of that gap
for the operator $ P_\nu $ obtained in \cite[Theorem 2]{damped} which was partly based on 
some estimates from \cite{nozw}. The required $ H^{s_0 } $ regularity in \eqref{eq:mixing} 
can be obtained by using ellipticity of $ \nu \Delta $ and that, roughly, accounts for the factor of
$ \nu^{-K} $ in \eqref{eq:theo}. 

Our setting and motivation is very different than that of \cite{enhance} and of the papers cited above.
But it is curious to ask if an improved version of the optimal enhanced dissipation involving the 
mixing constant \eqref{eq:theo}  is also true in time depending settings considered by those authors.
In the setting closer to ours one can also ask the same question for kinetic Brownian motion
-- see Drouot \cite{dro}. (For the ``large $ \nu $ limit" in that setting see Ren--Tao \cite{ret}.)

\noindent
{\sc Acknowledgements.} We would like to thank Jonathan Mattingly for introducing us to 
the problem of enhanced dissipation during a talk in Berkeley.
Our work has been partially supported by the
Simons Targeted Grant Award No. 896630, ``Moir\'e Materials Magic" .

\section{Preliminaries}
In this section we recall basic preliminaries in semiclassical analysis
and some specialized facts concerning spectral analysis of contact Anosov flows.

\subsection{Semiclassical analysis}
We will use some basic properties of semiclassical pseudodifferential operators on 
manifolds. They can all be found in \cite[Chapter 14]{ez}. On the flat space, 
a function $ a \in S^m (   \mathbb R^n \times 
\mathbb R^n ) $, that is a smooth function 
satisfying $ \partial^\alpha_x \partial^\beta_\xi a ( x, \xi ) = 
\mathcal O_{\alpha, \beta } ( \langle \xi \rangle^{m - |\beta | } ) $, $ 
\langle \xi \rangle := ( 1 + |\xi|^2 )^{\frac12} $, is quantized as
\[  a^{\rm{w}} ( x, h D ) u = \frac1 { (2 \pi h)^n } \int a ( ( x+ y)/2 , \xi ) e^{ \frac i h \langle x - y, \xi \rangle 
} u ( y ) d y d \xi . \] 
This can be generalized to manifolds using local charts and invariance under change of
coordinates to obtain $ \Psi_h^m ( M ) $, semiclassical pseudodifferential operators of order 
$ m $ on $ M $. 
We also have a quantization map $ S^m ( T^* M ) \ni a \to  \Op_h ( a ) :
H^{ s + m } ( M ) \to H^s ( M ) $ and a surjective symbol map $ \Psi_h^m ( M ) \ni A \to \sigma ( A )
\in S^m / hS^{m-1} ( T^* M ) $, $ \sigma ( \Op_h ( a ) ) = [ a ] $. An example we will use here
is given by functions of the Laplacian: for $ \chi \in C^\infty_{\rm{c}} ( \mathbb R ) $, 
and the Laplace-Beltrami operator for a metric $ g $, 
\[    \chi ( - h^2 \Delta_g  ) \in \Psi^{-\infty }_h ( M ) := \bigcap_{ m } \Psi_h^m ( M ) , \ \ \ 
\sigma (  \chi ( - h^2 \Delta_g  ) ) = [ \chi ( |\xi|_g^2 ) ]. \]
Let $\|u\|_{H^s_h}:=\|\langle -h^2\Delta\rangle^{s/2}u\|_{L^2}$, then for $A\in \Psi^m_h(M)$, there exists a constant $C$ independent of $h$ such that
\begin{equation*}
    \|Au\|_{H^s_h}\leq C\|u\|_{H^{s+m}_h}.
\end{equation*}
In applications, the symbols and functions may depend on $h$, we say a function depending on $h$ is tempered if there exists $N$ such that $\|u\|_{H^{-N}_h}\leq Ch^{-N}$. All functions considered below will be assumed to be tempered.

We also recall the basic composition formula:
$   \sigma_h ( A \circ B ) = \sigma_h ( A ) \sigma_h ( B ) $, or in more practical terms,
for $ a \in S^{m_1} $ and $ b \in S^{m_2} $, 
\begin{equation}
\label{eq:prod}  \Op_h ( a ) \circ \Op_h ( b ) = \Op_h ( c ) + R, \ \ \ 
R \in h^\infty \Psi^{-\infty}_h ( M ) , \ \ \ c - a b \in h S^{ m_1 + m_2 - 1}. 
\end{equation}
Finally we recall a microlocal elliptic estimate: if
$A\in \Psi_h^m(M)$, $B=\Op_h(b)\in \Psi_h^0(M)$ for some $b\in S^0(T^*M)$, and $T^*M=U\cup V$ for two open sets $U,V$ which are conic when $|\xi|$ is large, such that 
\begin{align*}
    \left\{ \begin{array}{l} 
|\sigma(A)| \geq c\langle \xi\rangle^{m}, \ \ (x,\xi)\in U,\\
|\partial_x^\alpha\partial_\xi^\beta b(x,\xi)|=\mathcal{O}_{\alpha,\beta}(h^\infty\langle \xi\rangle^{-\infty}), \ \ (x,\xi)\in V.
\end{array} \right. 
\end{align*}
Then
\begin{equation}\label{e:ell-micro}
    \|Bu\|_{H_h^s}\leq C\|Au\|_{H_h^{s-m}}+\mathcal{O}(h^\infty)\|u\|_{H^{-N}_h}.
\end{equation}
Taking $b=1$ we get a global ellipticity result:
\begin{equation}
\label{eq:elli}
A \in \Psi_h^m ( M) , \ \ | \sigma ( A ) | \geq c \langle \xi \rangle ^m \ 
\Longrightarrow 
\left\{ \begin{array}{l} 
\exists \, h_0 , B \in \Psi_h^{-m} ( M )  \ \forall \, 0 < h < h_0, \, s \in \mathbb R, \\
A \circ B = I_{ H^s ( M )},  \ \   B \circ A = I_{H^s ( M )  }. 
\end{array} \right. 
\end{equation}

\subsection{Hyperbolic dynamics}
A flow $\varphi^t=\exp(tX)$ on $M$ is Anosov if there is a continuous splitting of the tangent bundle $TM$:
\begin{equation*}
    T_xM=E_0(x)\oplus E_s(x)\oplus E_u(x)
\end{equation*}
where $E_0(x)=\mathbb{R}X(x)$ is the flow direction and there exist $C,\theta>0$ such that 
\begin{equation*}
    |d\varphi^t(x)v|_{\varphi^t(x)}\leq Ce^{-\theta|t|}|v|_x,\quad
\left\{ \begin{array}{l}
t>0,\ v\in E_s(x),\\
t<0,\ v\in E_u(x).
\end{array} \right. 
\end{equation*}
For an Anosov flow, we define the \emph{averaged Lyapunov exponent} as the largest constant $\gamma_0>0$ so that for any $\delta>0$,
\begin{equation}\label{e:lyapunov}
    |\det(d\varphi^{-t})|_{E_u(x)}|\leq C_\delta e^{-(\gamma_0-\delta)t},\quad t\geq 0.
\end{equation}
This exponent appears in the essential spectral gap (see \cite[Theorem 2]{damped}).

Now we switch to the spectral theory point of view.
We first observe that self-adjointness of $ i X$ on $ L^2 ( M ) $ (defined using the contact volume form)
implies that $ ( P_0 - \lambda)^{-1} : L^2 (M )  \to L^2 ( M ) $ for $ \Re \lambda < 0 $. 
We recall from \cite{FS} and \cite{zeta} that for any $ \gamma > 0 $ there exists $ s_0 >0  $ such that
this resolvent continues meromorphically to $ \Re \lambda < \gamma $ as an 
operator from $ H^{s_0 } ( M ) \to H^{-s_0 } ( M ) $. Abusing the notation slightly (one needs finer 
spaces to consider the actual inverse), 
\[ ( P_0 - \lambda)^{-1} :  H^{s_0 } ( M ) \to H^{-s_0 } ( M ) , \ \ \  \Re \lambda < \gamma , \]
with poles of finite rank. (For a precise discussion of the dependence of $ s_0 $ on $ \gamma $
see the recent paper by Dyatlov \cite{D23}.) By \cite[Theorem 4]{nozw} (see also \cite{liver} and \cite{Ts}), $ ( P_0 - \lambda )^{-1} $ has only finitely many poles in 
$ \Re \lambda < \gamma_0/2 -\delta $ for any $\delta>0$.  Hence we can define 
\begin{equation}
\label{eq:beta}
  \beta_0 := \sup \{ \beta  : \ \text{$ \lambda \mapsto ( P_0 - \lambda)^{-1}
$ is holomorphic for $ \Re \lambda < \beta $, $ \lambda \neq 0 $} \} . 
\end{equation}
We then have \eqref{eq:mixing} for $ \beta_1 < \beta_0 $ (typically for $ \beta_1= \beta_0 $).

For general Anosov flows \cite[Theorem 1]{damped} showed that in any compact subset
of $ \mathbb C $, the eigenvalues of $ P_\nu $ (defined on $ L^2 ( M ) $ with the domain 
$ H^2 (M ) $) converge to the poles of the meromorphic continuation of 
$ ( P_0 - \lambda)^{-1} $ that is to Pollicott--Ruelle resonances. In the case of contact
flows, \cite[Theorem 2]{damped} showed the uniformity of the spectral gap as $ \nu \to 0 $
with a quantitative resolvent estimate. Here we need a slight improvement:
\begin{equation}
\label{e:gap_res}
\begin{gathered}
\text{For any $ \beta < \min\{\beta_0,\gamma_0/2\}$ and $ \varepsilon > 0 $ there exist $ 
C,s_0, K_0 , \nu_0 >0$ such that } \\
         \|(P_\nu-\lambda)^{-1}\|_{H^{s_0}\to H^{-s_0}}\leq C\langle \lambda\rangle^{K_0}, \ \
   -1\leq \Re \lambda \leq \beta,  \ |\lambda | \geq \varepsilon, \ \  0 \leq \nu < \nu_0 . 
    \end{gathered}
          \end{equation}
 \begin{proof}[Proof of \eqref{e:gap_res}]
     The case $|\Im \lambda|> R$ follows from \cite[Theorem 2]{damped} and
      we only need to consider the region $|\Im \lambda|\leq R$. By \cite[Theorem 1]{damped}, for 
       $\beta < \beta_0 $ and $ \nu $ small enough, 
     $ \Spec ( P_\nu ) \cap \{ \Re \lambda \leq \beta\} = \{ 0 \} $, and the region $\{-1\leq \Re\lambda\leq \beta\}\setminus B_\epsilon(0)$ is uniformly away from $\Spec(P_\nu)$. In order to get the resolvent bound, we recall from \cite[Proposition 4.3]{damped} that there is $Q\in\Psi_h^{-\infty}$ such that for $0<\nu\leq  h/C_0\ll 1$ ($h$ will be fixed) we have
     \begin{equation*}
         \|(P_\nu+h^{-1}Q-\lambda)^{-1}\|_{H^{s_0}\to H^{-s_0}}\leq Ch^{-K_1}.
     \end{equation*}
     We also note that 
     \begin{equation}
     \label{eq:Q20}
         (P_\nu-\lambda)^{-1}=(P_\nu+h^{-1}Q-\lambda)^{-1}(I-K_\nu(\lambda))^{-1}
     \end{equation}
     where $K_\nu(\lambda)=h^{-1}Q(P_\nu+h^{-1}Q-\lambda)^{-1}$. By \cite[Lemma 5.1]{damped}, the map $\nu\mapsto K_\nu(\lambda)$ is continuous in $H^{s_0}\to H^{s_0}$ norm. So for $\nu$ sufficiently small, we have
  \begin{equation*}
      \|(I-K_\nu(\lambda))^{-1}-(I-K_0(\lambda))^{-1}\|_{H^{s_0}\to H^{s_0}}<1.
  \end{equation*}
  Since we are away from the poles of $ ( P_0 - \lambda)^{-1} $, $\|(I-K_0(\lambda))^{-1}\|_{H^{s_0}\to H^{s_0}}$ is uniformly bounded (see \eqref{eq:Q20} and note that $ h$ is small but fixed), it follows
  that  $\|(I-K_\nu(\lambda))^{-1}\|_{H^{s_0}\to H^{s_0}}$ is also uniformly bounded and hence
$    \|(P_\nu-\lambda)^{-1}\|_{H^{s_0}\to H^{-s_0}}\leq C$. 
 \end{proof}

\section{Proof of theorem}
%\begin{proof}[Proof of Theorem \ref{t:main theorem}]

    The proof follows from a contour deformation argument. We first recall that $ P_\nu $,
    with the domain given by $ H^2 ( M ) $, generates a strongly continuous semigroup
    $ e^{ - t P_\nu }$. This follows from the Hille--Yosida Theorem \cite[Theorem X.47a]{reesi}
    applied to $ P_\nu + C_0 \nu $, for a sufficiently large constant $ C_0 > 0 $. (If we choose 
    the metric so that $ \Delta $ is self-adjoint with respect to the contact volume form  
    we can take $ C_0 = 0 $). Also, for $ t > 0$,
    \begin{equation}
    \label{eq:cont1}
        e^{-tP_\nu}=\frac{1}{2\pi i}\int_{\Re \lambda=-1}(P_\nu-\lambda)^{-1}e^{-\lambda t}d\lambda, 
    \end{equation}
in the sense of distributions in $ t $: for $ \chi \in C^\infty_{ \rm{c}} ( ( 0 , \infty )) $, 
\begin{equation}
\label{eq:cont}
\int_0^\infty e^{-t P_\nu} \chi ( t ) dt =   \frac{1}{2\pi } \int_{\mathbb
R }(P_\nu + 1 - i \tau )^{-1} \widehat \chi ( \tau + i ) d\tau . 
\end{equation}

\begin{figure}
\includegraphics[width=11.5cm]{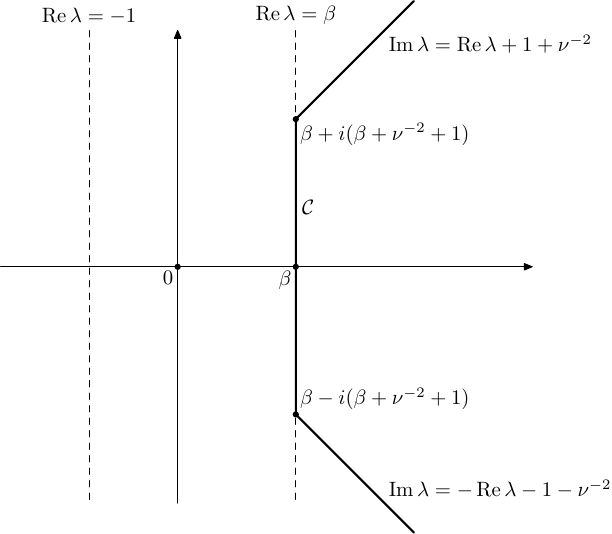}
\caption{\label{f:2} The contour of deformation used in the proof of the main theorem.}
\end{figure}

 To obtain \eqref{eq:theo}, we will deform the contour in \eqref{eq:cont1} to a contour shown in 
 Figure \ref{f:2},
    \begin{equation*}
         {\mathcal C}:=\{\Re\lambda=\beta,|\Im \lambda|\leq 1+\nu^{-2}+\beta\}\cup \{|\Im \lambda|= 1+\nu^{-2}+\Re\lambda,\Re\lambda\geq \beta\}.
    \end{equation*}    

        By \eqref{e:gap_res}, it suffices to show there is no eigenvalue of $P_\nu$ in the region between $\{\Im\lambda=\beta\}$ and $ {\mathcal C}$, and that for some $ N_1 $, 
\begin{equation}\label{e:resolventbound}
        \|(P_\nu-\lambda)^{-1}\|_{L^2\to L^2}\leq C\nu^{-N_1},\quad \lambda\in  {\mathcal C}.
    \end{equation}
In fact, once \eqref{e:resolventbound} is established, we have
\begin{equation*}
    e^{-tP_\nu}u=\underline u+\frac{1}{2\pi i}\int_{\mathcal  C}(P_\nu-\lambda)^{-1}e^{-\lambda t}u\, d\lambda.
\end{equation*}
The contour deformation (for $t>0$) is justified using \eqref{eq:cont} so that 
we get rapid decay in $\lambda$ (see step 1 below for the resolvent bound in the region \eqref{e:far_region} and note that $ \widehat \chi ( - i \lambda ) 
= \mathcal O ( e^{ C (-\Re \lambda)_+ } ( 1 + |\lambda |)^{-N }) $, for any $ N$). The main theorem then follows form
\begin{equation*}
    \frac{1}{2\pi}\int_{ {\mathcal C}}\|(P_\nu-\lambda)^{-1}\|_{L^2\to L^2}e^{-t\Re\lambda }\,d |\lambda|\leq C\nu^{-N_1} (\nu^{-2}+1)e^{-\beta t},\quad t>1. 
\end{equation*}

\begin{proof}[Proof of \eqref{e:resolventbound}]
\noindent
\textbf{Step 1:} We first deal with estimates on 
$$ {\mathcal C}_{\rm far}:=\{|\Im \lambda|=1+\nu^{-2}+\Re\lambda,\ \Re\lambda\geq \beta\}. $$ For that consider $\lambda\in{\mathbb C}$ with
\begin{equation}\label{e:far_region}
    \Re\lambda\geq -1,\quad |\Im \lambda|\geq 1+\nu^{-2}+|\Re\lambda|.
\end{equation}
For $f\in L^2$ and $u=(P_\nu-\lambda)^{-1}f$, we have
\begin{equation*}
    \left(-\nu|\lambda|^{-1}\Delta+|\lambda|^{-1}X-{\lambda}{|\lambda|^{-1}}\right)u=|\lambda|^{-1}f.
\end{equation*}
Putting  $h=\nu^{1/2}|\lambda|^{-1/2}$, we have  
\begin{equation*}
  -\nu|\lambda|^{-1}\Delta+|\lambda|^{-1}X-{\lambda}{|\lambda|^{-1}}\in \Psi^2_h , 
\end{equation*}
with principal symbol $|\xi|^2-\lambda/|\lambda|+\mathcal{O}(h^{1/3})_{S^1}$. This symbol
 is elliptic (see \eqref{eq:elli}) since $|\Im \lambda|>|\Re\lambda|$. As long as $\nu$ is sufficiently small, we conclude from \eqref{eq:elli} (the bounds are uniform on $ L^2 $ and on semiclassical Sobolev spaces)
\begin{equation*}
    \|u\|_{L^2}\leq C\left\|\left(-\nu|\lambda|^{-1}\Delta+|\lambda|^{-1}X-{\lambda}{|\lambda|^{-1}}\right)u\right\|_{L^2}\leq C|\lambda|^{-1}\|f\|_{L^2}.
\end{equation*}
Hence there are not eigenvalue of $P_\nu$ in the region defined by \eqref{e:far_region}, and \eqref{e:resolventbound} holds with $N_1=-2$.
    
    \noindent\textbf{Step 2:} We now deal with the estimate on 
     $$ {\mathcal C}_{\rm near}:=\{\Re\lambda=\beta,|\Im \lambda|\leq 1+\nu^{-2}+\beta\}.$$
      From \eqref{e:gap_res} we obtain
    \begin{equation*}
        \|(P_\nu-\lambda)^{-1}\|_{H^{s_0}\to H^{-s_0}}\leq C\langle \lambda\rangle^{K_0}\leq C \nu^{-2K_0}, \ \ 
        \lambda \in \mathcal C_{\rm near} . 
    \end{equation*}
  Note $\Re\lambda\leq \beta$ and
  \begin{align*}
      \|u\|_{H^s}^2&\lesssim \Re\langle -\Delta u , u\rangle_{H^{s-1}}+\|u\|_{H^{s-1}}^2\\
      &\lesssim \nu^{-1}\Re\langle (-\nu\Delta u+X-\lambda)u , u\rangle_{H^{s-1}}+\nu^{-1}\|u\|_{H^{s-1}}^2,
  \end{align*}
    we have 
    \begin{equation}\label{e:ell}
        \|u\|_{H^s}\lesssim \nu^{-1}\|(P_\nu-\lambda)u\|_{H^{s-1}}+\nu^{-1}\|u\|_{H^{s-1}}.
    \end{equation}
    By iterating the estimate \eqref{e:ell} we get an $H^{s_0} \to L^2$ bound (since we do not care about precise constants
    we can assume that $ s_0 $ is a positive integer): 
\begin{equation}\label{e:HstoL2}
    \begin{aligned}
            \|u\|_{L^2} \lesssim \nu^{-s_0}\|(P_\nu-\lambda)u\|_{L^2}+\nu^{- s_0}\|u\|_{H^{-s_0}}\lesssim \nu^{-s_0-2K_0}\|(P_\nu-\lambda)u\|_{H^{s_0}}.
    \end{aligned}
    \end{equation} 
    
    In order to deduce an $L^2\to L^2$ bound from this, we consider high frequencies and low frequencies separately.
    %We decompose $f\in L^2$ into $f=f_1+f_2$, where
    %$f_1=\chi(\nu^{2}D)f$ and $\chi$ is a cutoff function near $\{\xi=0\}$.
    %Then for the low frequency part
    %\begin{equation*}
    %    \|(P_\nu-\lambda)^{-1}f_1\|_{L^2}\lesssim \nu^{-s_0-2K_0}\|f\|_{H^{s_0}}\lesssim \nu^{-3s_0-2K_0} \|f\|_{L^2}.
    %\end{equation*}
    Let $f\in L^2$ and $u=(P_\nu-\lambda)^{-1}f$, so that 
    \begin{equation*}
        (-\nu^4\Delta +\nu^3 X-\nu^3\lambda)u=\nu^3 f.
    \end{equation*}
    We put $h=\nu^2$, so that  $-\nu^4\Delta +\nu^3 X-\nu^3\lambda\in \Psi^2_h$ with principal symbol $|\xi|^2+\mathcal{O}(h^{1/2})_{S^1}$, which is elliptic outside the zero section, $\{\xi=0\}$.
    Let $\chi \in C^\infty_{\rm c } ( \mathbb R ; [0 , 1 ] ) $  be equal to $ 1 $ near $ 0$. 
   The elliptic estimate \eqref{e:ell-micro} gives 
    \begin{equation*}
        \|(1-\chi(-\nu^4 \Delta)) u\|_{L^2}\lesssim \|\nu^3 f\|_{L^2}+ {\mathcal O}(\nu^\infty)\|u\|_{L^2}.
    \end{equation*}
    For the low frequency part,  $\chi(- \nu^4 \Delta)u$,  we use \eqref{e:HstoL2}:
    \begin{equation*}
    \begin{split} 
        \|\chi(- \nu^4 \Delta)u\|_{L^2}& \lesssim \nu^{-  s_0-2K_0}\|(P_\nu-\lambda)\chi(- \nu^4 \Delta)u\|_{H^{s_0}}\\
        & \lesssim \nu^{-  s_0-2K_0}(\|\chi(- \nu^4 \Delta)f\|_{H^{s_0}}+\|[\chi(- \nu^4 \Delta),P_\nu]u\|_{H^{s_0}}).
        \end{split} 
    \end{equation*}
   Since we can use $ \Delta $ to obtain an equivalent norm on Sobolev spaces, 
   \begin{equation*}
       \|\chi(- \nu^4 \Delta)f\|_{H^{s_0}}\lesssim \nu^{-2s_0}\|f\|_{L^2}.
   \end{equation*}
   For the commutator term we observe that $[\chi(- \nu^4 \Delta),\nu^3 P_\nu]$ is microlocalized near $\supp \chi'$, that is away from the zero section. Hence, we can use the elliptic estimate \eqref{e:ell-micro} again: 
    \begin{equation*}
        \|[\chi(- \nu^4 \Delta),P_\nu]u\|_{H^{s_0}}\lesssim \nu^{-2s_0}\|f\|_{L^2}+{\mathcal O}(\nu^\infty)\|u\|_{L^2}.
    \end{equation*}
    Putting these estimates together we conclude
       \[  \|(P_\nu-\lambda)^{-1}\|_{L^2\to L^2}\lesssim \nu^{-3s_0-2K_0} , \]
    which shows \eqref{e:resolventbound} with $ N_1 =3 s_0+2K_0$.\qedhere
\end{proof}

\end{document}